\documentclass[12pt]{article}
\usepackage{amsmath}
\usepackage{amscd}
\usepackage{mathrsfs}
\usepackage{amsfonts}
\usepackage{amsmath,amssymb}
\baselineskip 5.5pt \lineskip 5.5pt \numberwithin{equation}{section}

\def\qed{\hfill$\Box$\par}

\def \Z{\hbox{$Z\hskip -5.2pt Z$}}

\def \C{\hbox{$C\hskip -5pt \vrule height 6pt depth 0pt \hskip 6pt$}}

\def\qed{\ \ \ifhmode\unskip\nobreak\fi\ifmmode\ifinner
         \else\hskip5pt\fi\fi
 \hbox{\hskip5pt\vrule width4pt height6pt depth1.5pt\hskip 1 pt}}
\def\a{\alpha}

\def\cl{\centerline}

\def\vs{\vspace*}

\def\LL{{\cal L}}

\def\C{\mathbb{C}}

\def\Z{\mathbb{Z}}

\def\vs{\vspace*}
\def\cl{\centerline}
\def\QED{\hfill$\Box$}
\def\W{\mathcal {W}}

\newtheorem{theo}{Theorem}[section]
\newtheorem{lemm}[theo]{Lemma}
\newtheorem{rema}[theo]{Remark}
\newtheorem{defi}[theo]{Definition}
\newtheorem{coro}[theo]{Corollary}
\newtheorem{prop}[theo]{Proposition}
\begin{document}

\vs{10pt} \cl{\large {\bf Low dimensional cohomology of Hom-Lie
algebras}} \cl{\large{\bf and $q$-deformed $W(2,2)$
algebra}\footnote{Corresponding author: lmyuan@mail.ustc.edu.cn}}

 \vs{8pt}

\cl{ Lamei Yuan$^{\,\ddag}$, Hong You$^{\,\ddag\,\dag}$ }
 \cl{\small{ $^{\ddag}$Academy of Fundamental and Interdisciplinary
 Sciences,}}
\cl{\small{Harbin Institute of Technology, Harbin 150080, China}}
\cl{\small{$^{\dag}$Department of Mathematics, Suzhou University,
Suzhou 200092, China}} \cl{\small E-mail: lmyuan@mail.ustc.edu.cn,
youhong@suda.edu.cn
 }\vs{6pt}

\noindent{{\bf Abstract.} This paper aims to study the low
dimensional cohomology of Hom-Lie algebras and $q$-deformed $W(2,2)$
algebra. We show that the $q$-deformed $W(2,2)$ algebra is a Hom-Lie
algebra. Also, we establish a one-to-one correspondence between the
equivalence classes of one dimensional central extensions of a
Hom-Lie algebra and its second cohomology group, leading us to
determine the second cohomology group of the $q$-deformed $W(2,2)$
algebra. In addition, we generalize some results of derivations of
finitely generated Lie algebras with values in graded modules to
Hom-Lie algebras. As application we compute all $\a^k$-derivations
and in particular the first cohomology group of the $q$-deformed
$W(2,2)$ algebra. \vs{5pt}

\noindent{\bf Key words:} Hom-Lie algebras, $q$-deformed $W(2,2)$
algebra, derivation, second cohomology group, first cohomology
group.}

\noindent{\bf Mathematics Subject Classification (2000):} 17A30,
17A60, 17B68, 17B70.
\parskip .001 truein\baselineskip 8pt \lineskip 8pt

\vs{6pt}
\par
\noindent{\bf 1. \ Introduction}
\setcounter{section}{1}\setcounter{theo}{0}\setcounter{equation}{0}
\vs{6pt}

 The notion of Hom-Lie algebras was initially introduced
 in \cite{HLS} motivated by examples of
deformed Lie algebras coming from twisted discretizations of vector
fields. In this paper we will follow the slightly more general
definition of Hom-Lie algebras given by Makhlouf and Silvestrov in
\cite{MS1}. Precisely, a Hom-Lie algebra is a triple
$(\LL,[\cdot,\cdot],\alpha)$ consisting of a vector space $\LL$, a
bilinear map $[\cdot,\cdot]:\LL\times \LL\rightarrow \LL$ and a
linear map $\alpha:\LL\rightarrow \LL$ such that
\begin{eqnarray*}
[x,y]&=&-[y,x], \ \ (\mbox{skew-symmetry})\\
\circlearrowleft_{x,y,z}[[x,y],\alpha(z)]&=&0,\ \ \
(\mbox{Hom-Jacobi identity})
\end{eqnarray*}
for all $x,y,z\in \LL$, and where the symble
$\circlearrowleft_{x,y,z}$ denotes summation over the cyclic
permutation on $x,y,z.$ One sees that the classical Lie algebras
recover from Hom-Lie algebras if the twisting map $\a$ is the
identity map. The Hom-Lie algebras were discussed intensively in
\cite{MS2,SH,Yau1,Yau2} while the graded cases were considered in
\cite{AM,LS3,Y1}. But the cohomology with values in graded
Hom-modules is not very clear. Therefore, one of the aims of the
present paper is to fill this gap.

 The $W(2,2)$
Lie algebra was introduced in \cite{ZD} for the study of
classification of vertex operator algebras generated by vectors of
weight $2$. It is an extension of the Virasoro algebra. In the
following we denote by $\W$ the centerless $W(2,2)$ Lie algebra,
which is an infinite dimensional Lie algebra generated by $L_n$ and
$M_n$ ($n\in\Z$) satisfying the following Lie brackets
\begin{eqnarray*}
\ \ \ [L_m,L_n]=(n-m)L_{m+n}, \ {[L_m,M_n]}=(n-m)M_{m+n},\
{[M_m,M_n]}=0, \ \mbox{for}\ m,n\in\Z.
\end{eqnarray*}
In \cite{Y2} we presented a realization of the centerless $W(2,2)$
Lie algebra $\W$ by using bosonic and fermionic oscillators. The
bosonic oscillator $a$ and its hermitian conjugate $a^+$ obey the
commutation relations:
\begin{eqnarray}\label{boson}
[a,a^+]=aa^+-a^+a=1,\ \ \ [1,a^+]=[1,a]=0.
\end{eqnarray}
It follows by induction on $n$ that
$$[a,(a^+)^n]=n(a^+)^{n-1},\ \mbox{for all} \ n\in\Z.$$
The fermionic oscillators $b$ and $b^+$ satisfy the anticommutators
\begin{eqnarray}\label{fermion}
\{b,b^+\}=bb^++b^+b=1,\ \ \ b^2=(b^+)^2=0.
\end{eqnarray}
Moreover, we set $[a,b]=[a,b^+]=[a^+,b]=[a^+,b^+]=0.$
\begin{lemm}
{\rm(\it \cite{Y2})}\label{lemm-real} With notations above. The
generators of the form
\begin{eqnarray}\label{real}
\ \ \ \ \ \ \ \ L_n\equiv (a^+)^{n+1}a,\ \ \ M_n\equiv
(a^+)^{n+1}b^+a,\ \ \mbox{for all}\ n\in\Z,
\end{eqnarray}
realize the centerless $W(2,2)$ Lie algebra $\W$ under the
commutator
$$[A,B]=AB-BA,\ \ \mbox{ for all} \ A,B\in\W.$$
\end{lemm}

Now fix a nonzero $q\in\C$ such that $q$ is not a root of unity. We
introduce the following notation
\begin{eqnarray*}\label{notation}
[A,B]_{(\alpha,\beta)}=\alpha AB-\beta BA,
\end{eqnarray*}
and the $q$-number
\begin{eqnarray*}\label{qnumber}
[n]_q=\frac{q^n-q^{-n}}{q-q^{-1}}.
\end{eqnarray*}
It is clear to see that $[-n]_q=-[n]_q.$ Furthermore, one can deduce
that
\begin{eqnarray}\label{q}
q^{n}[m]_q-q^{m}[n]_q=[m-n]_q, \ \ \ q^{-n}[m]_q+q^{m}[n]_q=[m+n]_q.
\end{eqnarray}
Then the generators $L_n$ and $M_n$ ($n\in\mathbb{Z}$) satisfy the
following $q$-brackets:
\begin{eqnarray*}
[L_n,L_m]_{(q^{n-m},\,q^{m-n})}&=&[m-n]_qL_{m+n},\label{qLE1}\\
{[L_n,M_m]_{(q^{n-m},\,q^{m-n})}}&=&[m-n]_qM_{m+n},\label{qLE2}\\
{[M_n,M_m]_{(q^{n-m},q^{m-n})}}&=&0,\label{qLE3}
\end{eqnarray*}
for all $m,n\in\Z$. We call this algebra the {\it $q$-deformed
$W(2,2)$ algebra}, which is the second object considered in this
paper. In the following we will denote $q$-deformed $W(2,2)$ algebra
by $\mathcal {W}_q$ and simply write the $q$-bracket as
$[\cdot,\cdot]_q$. In \cite{Y2} we determined quantum groups and one
dimensional central extensions of $\mathcal {W}_q$. In this paper,
we will study its low dimensional cohomology theory. That is the
second aim of this paper.

 Throughout this paper, $\C$
denotes the field of complex number and $\mathbb{Z}$ denotes the set
of all integers. All
 vector spaces and algebras are assumed to be over $\C$.

\vskip12pt

\noindent{\bf 2. \ Second cohomology
group}\setcounter{section}{2}\setcounter{theo}{0}\setcounter{equation}{0}
\vskip6pt

In this section, we first recall some basic definitions and in
particular central extension of Hom-Lie algebras. Then we establish
a one-to-one correspondence between the equivalence classes of one
dimensional central extensions of a Hom-Lie algebra and its second
cohomology group with coefficients in $\C$. As application we
determine the second cohomology group of the $q$-deformed $W(2,2)$
algebra which is considered as a Hom-Lie algebra.

In the sequel we will often simply write a Hom-Lie algebra
$(\LL,[\cdot,\cdot],\a)$ as $(\LL,\a)$. A Hom-Lie algebra $(\LL,\a)$
is said to be {\it multiplicative} if the twisting map $\a$ is an
endomorphism. Let $G$ be an abelian group. A Hom-Lie algebra
$(\LL,\a)$ is said to be {\it $G$-graded}, if its underlying vector
space is $G$-graded (i.e., $\LL=\oplus_{g\in G}\LL_g$) satisfying
$[\LL_g,\LL_h]\subseteq \LL_{g+h}$, and if $\a$ is an even map,
i.e., $\a(\LL_g)\subseteq \LL_g$, for all $g,h\in G$.

   The theory of central extensions of Hom-Lie algebras was studied
in \cite{HLS,LS1}. An {\it extension} of a Hom-Lie algebra
$(\LL,\zeta)$ by an abelian Hom-Lie algebra $({\mathfrak
a},\zeta_{\mathfrak a})$ is a commutative diagram with exact rows
\begin{equation*}
\begin{CD}
0@>>> {\mathfrak a} @>{\rm \iota}>{\rm }> {\hat \LL} @>{\rm pr}>> \LL @>>> 0\\
@. @V{\zeta_{\mathfrak a}}VV @V{\hat \zeta}VV @V{\zeta}VV\\
0@>>> {\mathfrak a} @>{\rm \iota}>{\rm }> {\hat \LL} @>{\rm pr}>>
\LL @>>> 0
\end{CD}
\end{equation*}
where $({\hat \LL},{\hat \zeta})$ is a Hom-Lie algebra. The
extension is {\it central} if
$$\iota({\mathfrak a})\subseteq Z({\hat \LL})=\{x\in {\hat \LL}\,|\, [x,{\hat \LL}]_{\hat \LL}=0\}.$$

In the following we focus on the central extension of $(\LL,\a)$ by
a one-dimensional center $\mathbb{C} c$, where ${ c}=\iota(1)$. Note
that the center $\mathbb{C} c$ can be considered as the
one-dimensional trivial Hom-Lie algebra with the identity map.

\begin{defi}\rm \label{def-cocy}
Let $(\LL,\a)$ be a Hom-Lie algebra. A bilinear map $\psi: \LL\times
\LL\rightarrow \C$ is called a $2$-cocycle on $\LL$ if the following
conditions are satisfied
\begin{eqnarray}
\label{def-cocy-1} \!\!\!\!\!\!\!\!\!\!\!\!&&
\psi(x,y)=-\psi(y,x), \\
\label{def-cocy-2} \!\!\!\!\!\!\!\!\!\!\!\!&&
\psi(\a(x),[y,z])+\psi(\a(y),[z,x])+\psi(\a(z),[x,y])=0,
\end{eqnarray}
for all $x,y,z\in L$.
\end{defi}

Now we have the following theorem:
\begin{theo}\label{theo-cocycle} Let $(\LL,\a)$ be a Hom-Lie
algebra and $\psi: \LL\times \LL\rightarrow \C$ be a bilinear map.
Define on the vector space ${\hat \LL}=\LL\oplus \C$ the following
bracket and linear map by
\begin{eqnarray}
[x+c,y+b]_{\hat\LL}&=&[x,y]_{\LL}+\psi (x,y),\\
{\hat \a}(x+c)&=&\a(x)+c,
\end{eqnarray}
for all $x,y\in\LL$ and $c,b\in\C$. Then $({\hat
\LL},[\cdot,\cdot]_{\hat\LL},{\hat \a})$ is a Hom-Lie algebra one
dimensional central extension of $(\LL,\a)$ if and only if $\psi$ is
a $2$-cocycle on $(\LL,\a)$. If, in addition, $(\LL,\a)$ is
multiplicative and $\psi$ satisfies $\psi(\a(x),\a(y))=\psi(x, y)$,
for all $x,y\in\LL$, then the Hom-Lie algebra $({\hat \LL},{\hat
\a})$ is also multiplicative.
\end{theo}
\noindent{\it Proof.~} Since $[\cdot,\cdot]_{\LL}$ is
skew-symmetric, the new bracket $[\cdot,\cdot]_{\hat\LL}$ is
skew-symmetric if and only if the map $\psi$ is skew-symmetric. For
any $x,y,z\in\LL$ and $a,b,c\in\C$, we have
\begin{eqnarray*}
[{\hat\a}(x+a),[y+b,z+c]_{\hat\LL}]_{\hat\LL}&=&[\a(x)+a,[y,z]_{\LL}+\psi(y,z)]_{\hat\LL}\\
&=&[\a(x),[y,z]_{\LL}]_{\LL}+\psi(\a(x),[y,z]_{\LL}).
\end{eqnarray*}
Consequently,
\begin{eqnarray*}
[{\hat\a}(x+a),[y+b,z+c]_{\hat\LL}]_{\hat\LL}+[{\hat\a}(y+b),[z+c,x+a]_{\hat\LL}]_{\hat\LL}+[{\hat\a}(z+c),[x+a,y+b]_{\hat\LL}]_{\hat\LL}=0
\end{eqnarray*}
if and only if
\begin{eqnarray*}
\psi(\a(x),[y,z]_{\LL})+\psi(\a(y),[z,x]_{\LL})+\psi(\a(z),[x,y]_{\LL})=0,
\end{eqnarray*}
which proves $({\hat \LL},[\cdot,\cdot]_{\hat\LL},{\hat \a})$ is a
Hom-Lie algebra if and only if $\psi$ is a $2$-cocycle on
$(\LL,\a)$.

If $(\LL,\a)$ is multiplicative, then we have
\begin{eqnarray*}
{\hat\a}([x+a,y+b]_{\hat\LL})&=&{\hat\a}([x,y]_\LL+\psi(x,y))\\
&=&\a([x,y]_\LL)+\psi(x,y)\\
&=&[\a(x),\a(y)]_\LL+\psi(x,y).
\end{eqnarray*}
On the other hand, we have
\begin{eqnarray*}
[{\hat\a}(x+a),{\hat\a}(y+b)]_{\hat\LL}&=&[\a(x)+a,\a(y)+b]_{\hat\LL}\\
&=&[\a(x),\a(y)]_\LL+\psi(\a(x),\a(y)).
\end{eqnarray*}
According to the hypothesis that $\psi(\a(x),\a(y))=\psi(x, y)$ for
all $x,y\in\LL$, we have
\begin{eqnarray*}
{\hat\a}([x+a,y+b]_{\hat\LL})=
[{\hat\a}(x+a),{\hat\a}(y+b)]_{\hat\LL},\ \mbox{for\ all}\
x,y\in\LL,\ a,b\in\C,
\end{eqnarray*}
which shows that $({\hat \LL},{\hat \a})$ is multiplicative.

 Finally, we define $\rm pr$ and $\iota$ as the natural projection and
inclusion respectively by
\begin{eqnarray*}
\!\!\!\!\!\!\!\!\!\!\!\!&&
{\rm pr}:\hat \LL\rightarrow \LL,\qquad {\rm pr}(x+a)=x; \\
\!\!\!\!\!\!\!\!\!\!\!\!&& \iota:\C\rightarrow \hat \LL,\qquad
\iota(a)=0+a.
\end{eqnarray*}
Then it is easy to show that $({\hat \LL},{\hat \a})$ is a
one-dimensional central extension of $(\LL,\a)$.\QED\vs{2mm}

Denote by $Z^2(\LL,\C)$ the vector space of all $2$-cocycles on a
Hom-Lie algebra $(\LL,\a)$. For any linear map $f
:\LL\rightarrow\C$, we can define a $2$-cocycle $\psi_f$ by
\begin{equation}
\label{def-f}\psi_f(x,y)=f([x,y]),\ \ \mbox{for\ all}\ x,y\in \LL.
\end{equation}
Such a $2$-cocycle is called a {\it $2$-coboundary} or a {\it
trivial $2$-cocycle} on $\LL$. Let $B^2(\LL,\C)$ denote the vector
space of all $2$-coboundaries on $\LL$. The quotient space
$$H^2(\LL,\C) = Z^2(\LL,\C)/B^2(\LL,\C)$$
is called the {\it second cohomology group} of $\LL$ with trivial
coefficients $\C$. A $2$-cocycle $\psi$ is said to be equivalent to
another $2$-cocycle $\phi$ if $\psi-\phi$ is trivial. For a
$2$-cocycle $\psi$, let $[\psi]$ be the equivalent class of $\psi$.
Then we have the following corollary:
\begin{coro}\label{coro}For any Hom-Lie algebra $(\LL,\a)$, there exists
a one-to-one correspondence between the equivalence classes of one
dimensional central extensions of $(\LL,\a)$ and its second
cohomology group $H^2(\LL,\C)$.
\end{coro}

In the following, we consider the $q$-deformed $W(2,2)$ algebra
$\W_q$. Note that $\mathcal {W}_q$ is not a Lie algebra, because the
classical Jacobi identity does not hold (but the antisymmetry is
true). By straightforward calculations, we have
\begin{eqnarray}
&&(q^l+q^{-l})[\,[\,L_m,L_n]_{(q^{m-n},\,q^{n-m})},L_l\,]_{(q^{m+n-l},\,q^{l-m-n})}+\mbox{cyclic\
permutations}=0,\ \ \ \ \ \ \ \ \label{qJI1}\\
&&(q^l+q^{-l})[\,[\,L_m,L_n]_{(q^{m-n},\,q^{n-m})},M_l\,]_{(q^{m+n-l},\,q^{l-m-n})}+\mbox{cyclic\
permutations}=0.\label{qJI2}
\end{eqnarray}
Define on $\mathcal {W}_q$ a linear map $\alpha$ by
\begin{eqnarray*}
\alpha(L_n)=(q^n+q^{-n})L_n,\ \  \alpha(M_n)=(q^n+q^{-n})M_n.
\end{eqnarray*}
Then, using the $q$-deformed Jacobi identities (\ref{qJI1}) and
(\ref{qJI2}), we obtain the following result.
\begin{theo} \label{thom} The triple ($\mathcal {W}_q, [\cdot,\cdot]_q, \alpha$)
forms a Hom-Lie algebra.
\end{theo}

In \cite{Y2} we provided a computation of one-dimensional central
extensions of $\mathcal {W}_q$. Hence, according to Corollary
\ref{coro}, we can determine the second cohomology group of the
$q$-deformed $W(2,2)$ algebra $\mathcal {W}_q$ as follows:
\begin{prop}
$H^2(\mathcal {W}_q,\mathbb{C})=\mathbb{C}\beta\oplus
\mathbb{C}\gamma$, where
\begin{eqnarray*}
\beta(L_m,L_n)&=&\delta_{m,-n}\frac{[m-1]_q[m]_q[m+1]_q}{[2]_q[3]_q\langle
m\rangle_q},\ \ \ \beta(L_m,M_n)=\beta(M_m,M_n)=0,\\
\gamma(L_m,M_n)&=&\delta_{m,-n}\frac{[m-1]_q[m]_q[m+1]_q}{[2]_q[3]_q\langle
m\rangle_q},\ \ \ \gamma(L_m,L_n)=\gamma(M_m,M_n)=0,
\end{eqnarray*}
and where $\langle m\rangle_q=q^m+q^{-m}$, for all $m,n\in\Z.$
\end{prop}

\vskip10pt \noindent{\bf 3. \ Derivations of Hom-Lie algebras and
$q$-deformed $W(2,2)$ Lie
algebra}\setcounter{section}{3}\setcounter{theo}{0}\setcounter{equation}{0}
\vskip6pt

This section is devoted to discuss derivations of graded Hom-Lie
algebras. We extend to Hom-Lie algebras some concepts and results of
derivations of finitely generated Lie algebras with values in graded
modules studied in \cite{Fa}. As application we compute all
$\a^k$-derivations and particularly the first cohomology group of
the $q$-deformed $W(2,2)$ algebra.
\begin{defi}\rm Let $(\LL,\a)$ be a Hom-Lie algebra. A
representation of $\LL$ is a triple $(V, \rho,\beta)$, where $V$ is
a $\mathbb{C}$-vector space, $\beta\in End(V)$ and
$\rho:\LL\rightarrow End(V)$ is a $\mathbb{C}$-linear map satisfying
\begin{eqnarray*}
\rho([x,y])\circ
\beta=\rho(\a(x))\circ\rho(y)-\rho(\a(y))\circ\rho(x),
\end{eqnarray*}
for all $x,y\in \LL.$ $V$ is also called a Hom-$\LL$-module, denoted
by $(V,\beta)$ for convenience.
\end{defi}
One recovers the definition of a representation in the case of Lie
algebras by setting $\a={\rm id}_\LL$ and $\beta={\rm id}_V$. For
any $x\in \LL,$ define ${\rm ad}: \LL\rightarrow End(\LL)$ by ${\rm
ad}_x(y)=[x,y]$ for all $y\in \LL$. Then $(\LL,{\rm ad},\a)$  is a
representation of $\LL$, which is called  the {\it adjoint
representation }of $\LL$.
\begin{defi}
Let $(V,\beta_V)$ and $(W,\beta_W)$ be two Hom-$\LL$-modules. A
linear map $f:V\rightarrow W $ is called a morphism of
Hom-$\LL$-modules if it satisfies
\begin{eqnarray*}
f\circ \beta_V&=&\beta_W\circ f,\\
f(x\cdot v)&=&x\cdot f(v),
\end{eqnarray*}
for all $x\in\LL,$ $v\in V.$
\end{defi}

Let $G$ be an abelian group, $(\LL=\oplus_{g\in
G}\LL_g,[\cdot,\cdot],\a)$ be a $G$-graded Hom-Lie algebra. An
Hom-$\LL$-module $V$ is said to be $G$-graded if $V=\oplus_{g\in
G}V_g$ and $\LL_gV_h\subseteq V_{g+h}$ for all $g,h\in G$. For any
nonnegative integer $k$, denote by $\a^k$ the $k$-times composition
of $\a$, i.e.,
\begin{eqnarray}
\a^k=\a\circ\a\circ\cdots\circ\a.
\end{eqnarray}
In particular, $\a^0={\rm id}$ and $\a^1=\a.$ Then we can define
$\a^k$-derivations of $\LL$ with values in its Hom-$\LL$-modules.
\begin{defi} \rm A linear map $D:\LL\rightarrow V$ is called
an $\a^k$-derivation if it satisfies
\begin{eqnarray*}
D\circ \a&=&\a\circ D,\\
D[x,y]&=&\a^k(x)\cdot D(y)-\a^k(y)\cdot D(x),
\end{eqnarray*}
for all $x,y\in \LL$
\end{defi}

We recover the definition of a derivation by setting $k=0$ in the
definition above. Hence, an $\a^0$-derivation is often simply called
a derivation in the present paper. We say that an $\a^k$-derivation
$D$ has degree $g$ (denoted by ${\rm deg}(D)=g$) if $D\neq 0$ and
$D(\LL_h)\subseteq V_{g+h}$ for any $h\in G$. Let $D$ be an
$\a^k$-derivation. If there exists $v\in V$ such that
$D(x)=\a^k(x)\cdot v$ for all $x\in \LL$, then $D$ is called an
inner $\a^k$-derivation. Denote by $Der_{\a^k}(\LL,V)$ and
$Inn_{\a^k}(\LL,V)$ the space of $\a^k$-derivations and the space of
inner $\a^k$-derivations, respectively. In particular, let
$Der(\LL,V)$ and $Inn(\LL,V)$ denote the space of derivations and
the space of inner derivations, respectively, and write
$Der(\LL,V)_g:=\{D\in Der(\LL,V)\,\big|\,{\rm deg}(D)=g\}\cup\{0\}$.
The first cohomology group of $\LL$ with coefficients in $V$ is
defined by
\begin{eqnarray}
H^1(\LL,V):=Der(\LL,V)/ Inn(\LL,V).
\end{eqnarray}

\begin{rema}{\rm The set $Der_{\a^k}(\LL,V)$ (resp. $Inn_{\a^k}(\LL,V)$) is
not close under map composition or commutator bracket. But the space
of all such $\a^k$-derivations $\oplus_{k\geq 0}Der_{\a^k}(\LL,V)$
(resp. $\oplus_{k\geq 0}Inn_{\a^k}(\LL,V)$) form an Lie algebra via
commutator bracket.
 }
\end{rema}

Now let $\LL$ be a $G$-graded Hom-Lie algebra which is finitely
generated. In the following we present two results, which can be
seen as Hom versions of those obtained in \cite{Fa}.
\begin{prop}  Let $V$ be a $G$-graded
Hom-$\LL$-module. For every $D\in Der(\LL,V)$, we have
\begin{eqnarray}\label{summable}
D=\mbox{$\sum_{g\in G }D_g$},
\end{eqnarray}
where $D_g\in Der(\LL,V)_g$ and where there are only finitely many
$D_g(u)\neq 0$ in the equation $D(u)=\sum_{g\in G}D_g(u)$, for any
$u\in \LL$.
\end{prop}
\noindent{\it Proof}. For any $g\in G$, define a homogeneous linear
map $D_g:\LL\rightarrow V$ as follows: for any $u\in \LL_h$ with
$h\in G$, write $D(u)=\sum_{p\in G}u_p$ with $u_p\in V$, then set
$D_g(u)=u_{g+h}$. Clearly, $D_g$ is well defined and $D_g\in
Der(\LL,V)_{g}$. Also, (\ref{summable}) is true. \QED
\begin{prop}\label{Hom-Version} Let $V$ be a $G$-graded
Hom-$\LL$-module such that
\begin{itemize}\item[\rm(a)] $H^1(\LL_0,V_g)=0$, for $g\in G\setminus\{0\}$.
\item[\rm(b)]${\rm Hom}_{L_0}(\LL_g,V_h)=0$, for $g\neq h$.
\end{itemize}
Then $$Der(\LL,V)=Der(\LL,V)_0+Inn(\LL,V).$$
\end{prop}
\noindent{\it Proof}. Let $D$ be a derivation from $\LL$ into its
Hom-$\LL$-module $V$. According to (\ref{summable}) we can decompose
$D$ into its homogeneous components $D=\sum_{g\in G}D_g$ with
$D_g\in Der(\LL,V)_g.$ Suppose that $g\neq 0$. Then $D_g|_{\LL_0}$
is a derivation from $\LL_0$ into the Hom-$\LL_0$-module $V_g$. By
virtue of (a), $D_g|_{\LL_0}$ is inner, i.e., there exists $v_g\in
V_g$ such that $D_g(u)=u\cdot v_g$ for all $u\in \LL_0.$ Consider
$\psi_g:\LL\rightarrow V$ defined by $\psi_g(x):=D_g(x)-x\cdot v_g$
for all $x\in \LL$. Then $\psi_g$ is a derivation of degree $g$
which vanishes on $\LL_0$. Hence $\psi_g$ is a morphism of
Hom-$\LL_0$-modules and condition (b) entails the vanishing of
$\psi_g$ on $\LL_h$ for every $h\in G$. Consequently, $D_g\in
Inn(\LL,V)$, which completes the proof. \QED\vs{2mm}

In the following we focus on the $q$-deformed W(2,2) algebra
$\mathcal {W}_q$ as a Hom-Lie algebra ($\mathcal {W}_q,
[\cdot,\cdot]_q, \alpha$) defined in Theorem \ref{thom}. Obviously,
$\mathcal {W}_q$ is $\mathbb{Z}$-graded by
$$\mathcal {W}_q=\oplus_{n\in\mathbb{Z}}\mathcal {W}_q^n,\ \mbox{where}\ \mathcal {W}_q^n=span_\mathbb{C}\{L_n,M_n\}.$$
Note that $\mathcal {M}:=span_\mathbb{C}\{M_n\}$ is an ideal of
$(\mathcal {W}_q,\a)$, or in other words, $\mathcal {M}$ is an
adjoint Hom-$\mathcal{W}_q$-module. In addition, $\mathcal {W}_q$ is
finitely generated by $\{L_1,L_{-1},M_1\}$. Let $D$ be an
$\a^k$-derivation of $\mathcal {W}_q$. For all $m,n\in\mathbb{Z}$,
we have
\begin{eqnarray}
(q^m+q^{-m})^k[D(L_n),L_m]_q+(q^n+q^{-n})^k[L_n,D(L_m)]_q&=&[m-n]_qD(L_{m+n}),\label{qj11}\\
(q^m+q^{-m})^k[D(L_n),M_m]_q+(q^n+q^{-n})^k[L_n,D(M_m)]_q&=&[m-n]_qD(M_{m+n}).\label{qj22}
\end{eqnarray}
Now we aim to determine all $\a^k$-derivation of $\mathcal {W}_q$.
First, we compute the ($\a^0$-)derivations of $\mathcal {W}_q$.
Denote by $Der(\mathcal {W}_q)$ and $Inn(\mathcal {W}_q)$ the set of
all derivations and the set of all inner derivations, respectively.
Let $Der(\mathcal {W}_q)_m$ be the set of derivations of degree $m$.
\begin{lemm}$H^1(\mathcal {W}_q^0,\mathcal {W}_q^n)=0$ for any nonzero integer $n$.
\end{lemm}
\noindent{\it Proof}. Note that $[L_0,X_0]_q=0$, for any $X_0\in
\mathcal {W}_q^0=span\{L_0,M_0\}$. Let $D$ be any element in
$Der(\mathcal {W}_q^0,\mathcal {W}_q^n)$. Then it follows $D(L_0)\in
\mathcal {W}_q^n.$ Applying $D$ to $[L_0,X_0]_q=0$, we have $[n]_q
D(X_0)=[L_0,D(X_0)]=[X_0,D(L_0)]$. Consequently, $D(X_0)=[X_0,v]$
with $v=\frac{1}{[n]_q}D(L_0)$ in $\mathcal {W}_q^n$. In other
words, $D$ is an inner derivation from $\mathcal {W}_q^0$ into its
adjoint module $\mathcal {W}_q^n$.\QED
\begin{lemm}${\rm Hom}_{\mathcal {W}_q^0}(\mathcal {W}_q^m,\mathcal {W}_q^n)=0$ for $m\neq n$.
\end{lemm}
\noindent{\it Proof}. Let $f\in{\rm Hom}_{\mathcal {W}_q^0}(\mathcal
{W}_q^m,\mathcal {W}_q^n)$ with $m\neq n$. For any $X_m\in\mathcal
{W}_q^m$, we have
$$(q^n+q^{-n})\big(f(X_m)\big)=\a\big(f(X_m)\big)=f\big(\a(X_m)\big)=(q^m+q^{-m})f(X_m),$$
leading to $f(X_m)=0$, since $m\neq n$. Hence, $f=0$. \QED\vs{2mm}

Now according to Proposition \ref{Hom-Version}, we have the
following result:
\begin{prop}\label{p22} $Der(\mathcal {W}_q)=Der(\mathcal {W}_q)_0+Inn(\mathcal
{W}_q)$.
\end{prop}

Thanks to Proposition \ref{p22}, the study of $Der(\mathcal {W}_q)$
reduces to that of its constitute of degree zero. Let $D$ be an
element of $Der(\mathcal {W}_q)_0$. For any integer $n$, assume that
\begin{eqnarray}
D(L_n)=a_n L_n+b_n M_n, \ \ D(M_n)=c_n L_n+d_n M_n,
\end{eqnarray}
where the coefficients are complex numbers. Applying $D$ to
$[L_0,L_n]_q=[n]_qL_n$, one can obtain $D(L_0)=0,$ i.e.,
$a_0=b_0=0$.
 Using (\ref{qj11}), we have
\begin{eqnarray*}
a_{m+n}=a_m+a_n,\ \ b_{m+n}=b_m+b_n,\ \mbox{for}\ m\neq n.
\end{eqnarray*}
Let $m=-n$. Then we have
\begin{eqnarray*}
 a_{-m}=-a_m,\ \ b_{-m}=-b_m, \ \mbox{for}\ m\neq 0.
\end{eqnarray*}
Furthermore, we have
\begin{eqnarray*}
 a_m=ma_1,\ \ b_m=mb_1, \ \ \mbox{for all}\ m\in\Z.
\end{eqnarray*}
Similarly, using (\ref{qj22}) we have
\begin{eqnarray*}
 c_{m+n}=c_n,\ \ d_{m+n}=a_m+d_n, \ \ \mbox{for}\ m\neq n,
\end{eqnarray*}
from which it follows
\begin{eqnarray*}
 c_{m}=c_0,\ \ d_{m}=ma_1+d_0, \ \ \mbox{for\ all}\ m\in\mathbb{Z}.
\end{eqnarray*}
Applying $D$ to $[M_1,M_0]=0$, we have $c_0=0$. It follows that
$c_m=0$ for all $m\in\Z$. Hence, there exist $a, b, d\in\C$ such
that
\begin{eqnarray}\label{D}
 D(L_{n})=n(aL_n+bM_n),\ \ D(M_{n})=(na+d)M_n, \ \ \mbox{for\ all}\
 n\in\mathbb{Z}.
\end{eqnarray}
From the discussions above we obtain the following result:
\begin{prop}All the derivations of $\mathcal {W}_q$ is
$$Der(\mathcal {W}_q)=span_{\C}\{D\}\oplus Inn(\mathcal
{W}_q),$$ where $D$ is defined by (\ref{D}).
\end{prop}
\begin{coro} The first cohomology group of $\,\mathcal{W}_q$ with
values in its adjoint module is one-dimensional.
\end{coro}

Next we compute the $\a^1$-derivations of $\mathcal {W}_q$. Let $D$
be an $\a^1$-derivation of degree $s$. Assume that
\begin{eqnarray*}
D(L_n)=a_{s,n} L_{n+s}+b_{s,n} M_{n+s}, \ \ D(M_n)=c_{s,n}
L_{n+s}+d_{s,n} M_{n+s},
\end{eqnarray*}
where the coefficients are complex numbers. Then from equation
(\ref{qj11}) we obtain
\begin{eqnarray}
[m-n]_qa_{s,m+n}&=&(q^m+q^{-m})[m-s-n]_qa_{s,n}+(q^n+q^{-n})[m+s-n]_qa_{s,m},\label{qj33}\\
{[m-n]_qb_{s,m+n}}&=&(q^m+q^{-m})[m-s-n]_qb_{s,n}+(q^n+q^{-n})[m+s-n]_qb_{s,m}.\label{qj44}
\end{eqnarray}
We first consider the case of $s\neq 0$. Taking $m=0$ in
(\ref{qj33}), we have
\begin{eqnarray*}
(2[s+n]_q-[n]_q)a_{s,n}=(q^n+q^{-n})[s-n]_qa_{s,0}.
\end{eqnarray*}
Furthermore,
\begin{eqnarray}
a_{s,n}=\frac{(q^n+q^{-n})[s-n]_q}{(2[s+n]_q-[n]_q)}a_{s,0}.\label{88}
\end{eqnarray}
Plugging (\ref{88}) into (\ref{qj33}), we have
\begin{eqnarray*}
\frac{(q^{m+n}+q^{-m-n})[m-n]_q[s-m-n]_q}{2[s+m+n]_q-[m+n]_q}a_{s,0}&=&
\frac{(q^m+q^{-m})(q^n+q^{-n})[m-s-n]_q[s-n]_q}{2[s+n]_q-[n]_q}a_{s,0}\\
&+&
\frac{(q^m+q^{-m})(q^n+q^{-n})[m+s-n]_q[s-m]_q}{2[s+m]_q-[m]_q}a_{s,0}.
\end{eqnarray*}
Let $m=s$ in the equation above, we have
\begin{eqnarray}
\frac{(q^{s+n}+q^{-s-n})[s-n]_q[-n]_q}{2[2s+n]_q-[s+n]_q}a_{s,0}=\frac{(q^s+q^{-s})(q^n+q^{-n})[s-n]_q[-n]_q}{2[s+n]_q-[n]_q}a_{s,0}\label{ss}.
\end{eqnarray}
Then taking $n=-s$ in (\ref{ss}), we get
\begin{eqnarray*}
\frac{[2s]_q[s]_q}{[s]_q}a_{s,0}=\frac{(q^s+q^{-s})^2[2s]_q[s]_q}{[s]_q}a_{s,0}.
\end{eqnarray*}
It follows $a_{s,0}=0$ since $s\neq 0$. Then we have $a_{s,n}=0$ for
$n\in \mathbb{Z}$ and $s\neq 0$ by (\ref{88}). Similarly, from
(\ref{qj44}) we can deduce that $b_{s,n}=0$ for $s\neq 0$ and
$n\in\Z$.

In the case of $s=0$, we simply write $a_{0,n}$ as $a_n$. Then it
can be deduced from (\ref{qj33}) that
\begin{eqnarray}
a_{m+n}=(q^m+q^{-m})a_n+(q^n+q^{-n})a_m,\ \mbox{for}\ m\neq
n.\label{22}
\end{eqnarray}
Let $m=0$ in (\ref{22}), we have
\begin{eqnarray}
a_n=-(q^n+q^{-n})a_0,\label{22_1}
\end{eqnarray}
which implies $a_n=a_{-n}$ for $n>0$. Taking $m=-n$ in (\ref{22}),
we have
\begin{eqnarray}
a_0=(q^n+q^{-n})(a_n+a_{-n}).\label{22_2}
\end{eqnarray}
Substituting (\ref{22_1}) into (\ref{22_2}), we have $a_0=0$ and
$a_n=0$ for all $n\in\Z$. Similarly, we can deduce that $b_{0,m}=0$
for all $m\in\Z$ by using (\ref{qj44}) where $s=0$. Hence, we have
proved that
\begin{eqnarray*}
a_{s,m}=b_{s,m}=0, \ \mbox{for all}\ m,s\in\Z,
\end{eqnarray*}
or, in other words, we get $D(L_m)=0$ for $m\in\Z.$

It remains to determine $D(M_n)$ for all $n\in\Z.$ Using $D(L_n)=0$,
we can deduce from (\ref{qj22}) that
\begin{eqnarray}
[m-n]_qc_{s,m+n}&=&(q^n+q^{-n})[m+s-n]_qc_{s,m},\label{qj55}\\
{[m-n]_qd_{s,m+n}}&=&(q^n+q^{-n})[m+s-n]_qd_{s,m}.\label{qj66}
\end{eqnarray}
Let $m=0$ in (\ref{qj55}) and (\ref{qj66}), respectively. We have
\begin{eqnarray}
-[n]_q c_{s,n}&=&(q^n+q^{-n})[s-n]_qc_{s,0},\label{qj44_1}\\
-[n]_q d_{s,n}&=&(q^n+q^{-n})[s-n]_qd_{s,0}.\label{qj44_2}
\end{eqnarray}
Taking $n=0$ in (\ref{qj55}) and (\ref{qj66}), respectively, one has
\begin{eqnarray}
[m]_qc_{s,m}&=&2[m+s]_qc_{s,m}.\label{33_1}\\
{[m]_qd_{s,m}}&=&2[m+s]_qd_{s,m}.\label{33_2}
\end{eqnarray}
Taking $m=0$ in (\ref{33_1}) and (\ref{33_2}), respectively, we have
$c_{s,0}=d_{s,0}=0$ for $s\neq 0$. Then it follows from
(\ref{qj44_1}) (resp. (\ref{qj44_2})) that $c_{s,n}=0$ (resp.
$d_{s,n}=0$) for $n\in\Z$ and $s\neq 0.$

If $s=0$, then it follows from (\ref{qj55}) that
\begin{eqnarray}
[m-n]_qc_{0,m+n}&=&(q^n+q^{-n})[m-n]_qc_{0,m}.\label{77}
\end{eqnarray}
Let $n=0$ in (\ref{77}), we have $[m]_qc_{0,m}=2[m]_qc_{0,m}.$ It
follows that $c_{0,m}=0$ for $m\neq 0.$ Taking $n=-m$ in (\ref{77}),
we have $[2m]_qc_{0,0}=(q^m+q^{-m})[2m]_qc_{0,m}$, leading us to
$c_{0,0}=0$. Similarly, we can deduce from (\ref{qj66}) that
$d_{0,m}=0$ for all $m\in\Z.$ Thereby, the following proposition is
proved.
\begin{prop} If $D$ is an $\a^1$-derivation of $\mathcal {W}_q$,
then $D=0$.
\end{prop}

With the similar discussions as above, we can compute
$\a^k$-derivations of $\mathcal{W}_q$ for $k>1$ and thus we have all
the $\a^k$-derivations of $\mathcal{W}_q$ for $k>0$ determined.
\begin{prop}For $k>0$, all the $\a^k$-derivations of $\mathcal{W}_q$ are zero.
\end{prop}


\end{document}